\documentclass[11pt,twoside]{article}

\usepackage{graphicx}
\usepackage{multirow}
\usepackage{fancyhdr}
\usepackage{array}
\usepackage{amsfonts}
\usepackage{amsmath}
\usepackage{amssymb}
\usepackage{amsthm}

\usepackage{subfigure}

\usepackage{graphics} 
\usepackage{epsfig} 
\usepackage{latexsym}

\usepackage[textwidth=16cm,textheight=22cm,top=3cm,bottom=3cm]{geometry}
\usepackage{fancyhdr}

\newtheorem{theorem}{Theorem}
\newtheorem{lemma}{Lemma}

\newtheorem{remark}{Remark}
\theoremstyle{definition}
\newtheorem{definition}{Definition}

\usepackage{color}
\newenvironment{redtext}{}{\ignorespacesafterend}

\newcommand{\drop}[1]{}
\newcommand{\no}{\noindent}
\newcommand{\fer}[1]{(\ref{#1})}
\newcommand{\qtext}[1]{\quad\text{#1}}

\newcommand{\bx}{\mathbf{x}}
\newcommand{\by}{\mathbf{y}}

\newcommand{\cL}{\mathcal{L}}

\newcommand{\eps}{\varepsilon}
\newcommand{\vfi}{\varphi}

\newcommand{\grad}{\nabla}

\newcommand{\p}{\partial}

\newcommand{\N}{\mathbb{N}}
\newcommand{\R}{\mathbb{R}}

\def\O{\Omega}

\newcommand{\e}{{\text{e}}}

\newcommand{\abs}[1]{| #1 |}
\newcommand{\nor}[1]{\| #1 \|}

\title{ Well-posedness of an evolution problem with nonlocal diffusion   
\thanks{Supported by Spanish MCI Project MTM2017-87162-P.}}

\author{Gonzalo Galiano  \thanks{Dpt. of Mathematics, Universidad de Oviedo,
 c/ Calvo Sotelo, 33007-Oviedo, Spain ({\tt galiano@uniovi.es)}}}

\date{}

\begin{document}

\thispagestyle{plain}

\maketitle

\begin{abstract}
  We prove the well-posedness of a general evolution reaction-nonlocal diffusion problem under two sets of assumptions. In the first set, the main hypothesis is the Lipschitz continuity  of the range kernel and the bounded variation of the spatial kernel and the initial datum. In the second set of assumptions, we relax the Lipschitz continuity of the range kernel to H\"older  continuity, and assume  monotonic behavior. In this case, the spatial kernel and the initial data can be just integrable functions. The main applications of this model are related to the fields of Image Processing and Population Dynamics.

\no\emph{Keywords: }
Nonlocal diffusion, existence, uniqueness, $p-$Laplacian, bilateral filter.


\end{abstract}

\section{Introduction}

In this article, we study the well-posedness of a general class of evolution reaction-nonlocal diffusion problems expressed in the following form.  
Let $T>0$ and $\O\subset\R^d$ $(d\geq 1)$ be an open and bounded set with Lipschitz continuous boundary.  Find $u:[0,T]\times\overline{\O}\to\R$ such that 
\begin{align}
&  \p_t u(t,\bx)  =  \int_\O J(\bx-\by) A (t,\bx,\by,u(t,\by)-u(t,\bx)) d\by +f(t,\bx,u(t,\bx)), \label{eq.eq}\\
&   u(0,\bx)=  u_0(\bx),  \label{eq.id}
\end{align}
 for $(t,\bx)\in Q_T=(0,T)\times\O$, and for some $u_0:\O\to\R$. 

The main examples we have on mind are connected to the fields of Population Dynamics and of Image Processing. In the first case, choosing for instance $A(t,\bx,\by,s)=s$, we describe the balance   of population coming in and leaving from $\bx$, as 
\begin{align*}
\int_\O J(\bx-\by) u(t,\by) d\by -u(t,\bx),
\end{align*}
where the convolution kernel $J\geq0$, with $\int J=1$, determines the size and the shape of the  influencing neighborhood of $\bx$. 
In absence of a reaction term, the resulting equation is a nonlocal diffusion variant of the heat equation, usually written as  
\begin{align*}
 \p_t u(t,\bx)  = \int_\O J(\bx-\by) \big(u(t,\by)-u(t,\bx)\big) d\by .
\end{align*}
In this context, the nonlocal $p-$Laplacian diffusion operator, corresponding to $A(t,\bx,\by,s)=\abs{s}^{p-2}s$, for $p\in[1,\infty]$, is also a well known example, leading to the equation 
\begin{align*}
 \p_t u(t,\bx)  = \int_\O J(\bx-\by) \left|u(t,\by)-u(t,\bx)\right|^{p-2}\big(u(t,\by)-u(t,\bx)\big) d\by.
\end{align*}
These two examples correspond to a choice for which $A$ is a non-decreasing function of $s$. These kind of problems have been studied at great length by Andreu et al. in a series of works, see the monograph \cite{Andreu2010}. Their results, strongly dependent on the monotonicity of $A$, include the well-posedness as well as  properties such as the stability with respect to the initial data or the convergence of related rescaled nonlocal problems to their corresponding local versions. It is worth mentioning that problems of the type \fer{eq.eq} related to monotone functions, $A$, can be seen as gradient descents of convex energies. For instance, for the $p-$Laplacian, the nonlocal energy is given by
\begin{align*}
J_p (u) = \frac{1}{p} \int_\O\int_\O J(\bx-\by) \left|u(t,\by)-u(t,\bx)\right|^{p}.
\end{align*}

In the examples arising in Image Processing, the monotonicity of $A$ is not the rule. A very useful denoising filter, the bilateral filter \cite{Smith1997, Tomasi1998, Buades2005, Buades2006}, which provides results similar to the Perona-Malik equation \cite{Perona1990, Guidotti2014} or to the Total Variation restoration filter \cite{Rudin1992, {Chambolle2010}}, takes the form 
\begin{align*}
Bu(\bx) = \frac{1}{C(\bx)} \int_\O \exp\Big(-\frac{\abs{\bx-\by}^2}{\rho^2}\Big) \exp\Big(-\frac{\abs{u(\bx)-u(\by)}^2}{h^2}\Big)u(\by) d\by,
\end{align*}
 where $u$ is the image to be filtered, $\rho$ and $h$ are constants modulating the sizes of the space and range neighborhoods where the filtering process takes place, and $C$ is the normalizing factor 
\begin{align*}
C(\bx) = \int_\O \exp\Big(-\frac{\abs{\bx-\by}^2}{\rho^2}\Big) \exp\Big(-\frac{\abs{u(\bx)-u(\by)}^2}{h^2}\Big) d\by.
\end{align*}
Neighborhood filters like $B$ may also be derived from variational principles \cite{Kindermann2005}, being their correspondent gradient descent approximations  given by nonlocal equations of the type \fer{eq.eq}. Indeed, defining 
\begin{align}
\label{def.bil}
J(\bx)= \exp\Big({-\frac{\abs{\bx}^2}{\rho^2}}\Big), \quad A(t,\bx,\by,s) =  \exp\Big({-\frac{s^2}{h^2}}\Big) ,
\end{align}
we have that \fer{eq.eq} is the gradient descent associated to the nonconvex energy functional 
\begin{align*}
J_B (u) = \int_\O \int_\O \exp\Big(-\frac{\abs{\bx-\by}^2}{\rho^2}\Big)  \Big(1-  \exp\Big(-\frac{\abs{u(\bx)-u(\by)}^2}{h^2}\Big) \Big) d\bx d\by,
\end{align*}
for which the filter $Bu(\bx)$ is just a one step algorithm in the search direction.

From the definition of $A$ given in \fer{def.bil}, we readily see its lack of monotonicity. Thus, the approach followed by Andreu et al. may not be employed to show the well-posedness of the related gradient descent problem. 

Besides, there are other situations that we would like to cover for this kind of nonlocal diffusion problems which have been not treated, as far as we know, in the literature. 
One of them is allowing the convolution kernel, $J$, to be discontinuous. This is the case we encounter for  the Yaroslavsky filter \cite{Yaroslavsky1985}, with much faster numerical implementations than that of \fer{def.bil}, see \cite{Galiano2015b, Yang2015, Galiano2016}, which is given by 
\begin{align}
\label{def.yar}
J(\bx)= 1_{B_\rho(\bx)}(\by), \quad A(t,\bx,\by,s) =  \exp\Big({-\frac{s^2}{h^2}}\Big) ,
\end{align}
where $1_{B_\rho(\bx)}$ is the characteristic function of the ball $B_\rho(\bx)$.

Another situation we are interested in is that in which the power, $p$, of the $p-$Laplacian is not constant,  which finds applications in image restoration. Two important examples are the following:
\begin{enumerate}
\item The power $p$ depends on the nonlocal gradient of the solution, that is $p\equiv p(u(t,\by)-u(t,\bx))$, for some  even 
non-increasing smooth function $p:\R\to\R$. Typically, one takes $p(s)$ such that $p(0)=2$, leading to linear diffusion for small values of the nonlocal gradient of $u$, and $p(s)\to1$ as $s\to\infty$, implementing in this case the nonlocal Total Variation minimization in regions around sharp edges of $u$. Observe that the corresponding local version of this variational problem is the minimization of the functional introduced by Blomgren et al.  \cite{Blomgren1997},
\begin{align*}
J_G(u)=\int_\O \abs{\grad u(\bx)}^{p(\abs{\grad u(\bx)})-2}\grad u(\bx)d\bx.
\end{align*}
For the nonlocal version of $J_G$, the resulting gradient descent is given, in the form of equation \fer{eq.eq}, by 
\begin{align}
\label{eq:varplap}
 \p_t u(t,\bx)  = \int_\O J(\bx-\by) \left|u(t,\by)-u(t,\bx)\right|^{p(u(t,\by)-u(t,\bx))-2} \nonumber \\
 \times \big(u(t,\by)-u(t,\bx)\big) d\by,
\end{align}
see Remark~\ref{remark1}.

\item Since proving the existence of a minimum of $J_G$ is not evident, simpler functionals with $p\equiv p(\bx)$ have been introduced \cite{Chen2005, Guo2011}, with, for instance, $p(\bx)\equiv p(\grad u^\sigma_0(\bx))$, being $u^\sigma_0$ a regularization of the inital image.
As long as $u_0\in L^\infty(\O)\cap BV(\O)$, its nonlocal counterpart does not need of such regularization,    
and the resulting gradient descent takes the form 
\begin{align*}
 \p_t u(t,\bx)  = \int_\O J(\bx-\by) \left|u(t,\by)-u(t,\bx)\right|^{\tilde p(\bx,\by)-2}\big(u(t,\by)-u(t,\bx)\big) d\by.
\end{align*}
with $\tilde p(\bx,\by) \equiv p(u_0(\by)-u_0(\bx))$, and $p$ satisfying the properties  aforementioned.
\end{enumerate}

Turning to the results proved in this paper, the first theorem establishes the well posedness of problem \fer{eq.eq}-\fer{eq.id} under the main assumptions of Lipschitz continuity of $A(t,\bx,\by,s)$ as a function of $(\bx,\by,s)$, and of $L^\infty \cap BV$ regularity of the initial data. The latter is the usual regularity assumed for image representation, due to the convenience of allowing  discontinuities through level lines, which are the way in which  object edges are represented within the image. The Lipschitz continuity condition is satisfied, for instance, by the bilateral type filters and the nonlocal ${p(u(t,\by)-u(t,\bx))}$-Laplacian, see Remark~\ref{remark1}.

Although, in contrast to local diffusion, one of the relevant properties of nonlocal diffusion is its lack of a regularizing effect, we show that the regularity of the initial data and of the function $A$ are enough to 
use a compactness technique to deduce the existence of solutions.  The solution lies, therefore, in the same space, $L^\infty \cap BV$, that the initial data which is, in fact, a good property for image processing transformations. 

In our second result we relax the regularity of $A$ to just H\"older continuity with respect to $s$, but 
we additionally assume its monotonic behavior. This is the case of the nonlocal $p-$Laplacian for $p\in(1,\infty)$, extensively analyzed by Andreu et al. \cite{Andreu2010}. However, with our approach we are also able to deal with the $p(t,\bx,\by)-$Laplacian and with very general spatial kernels, $J$.  In addition, our result is based on a constructive technique that gives clues for the discrete approximation of such non-regular problems, see \cite{Karami2018} for a complementary approach.

\section{Assumptions and main results}  

Since $\O\subset\R^d$ is bounded, we have $\bx-\by \in B$ for all $\bx,\by\in\O$,
 for some open ball $B\subset\R^d$ centered at the origin. Thus, for $J$ defined on $\R^d$, we may 
 always replace it in \fer{eq.eq} by its restriction to $B$, $J|_B$. Abusing on notation, we write $J$ instead of  
$J|_B$ in the rest of the paper, and assume that $J$ is defined in a bounded domain.

We always assume, at least, the following hypothesis on the data.
  
\no \textbf{Assumptions (H)}
\begin{enumerate}
 \item The \emph{spatial kernel} $J\in L^1(B)$ is even and non-negative, with 
 \begin{align}
 \label{cond.j1}
 \int_{B} J(\bx)d \bx = 1.
 \end{align}

\item The \emph{range kernel} $A \in L^\infty((0,T)\times \O\times\O)\times C^{0,\alpha}_{loc}(\R)$, with $\alpha\in(0,1)$, satisfies:
\begin{align*}
A(\cdot,\cdot,\cdot,-s) = - A(\cdot,\cdot,\cdot,s),\quad  A(\cdot,\cdot,\cdot,s)s\geq0
\end{align*}
in $(0,T)\times\O\times\O$,  for all $s\in\R$, and 
\begin{align}
\label{ass:asym}
& A(\cdot,\bx,\by,\cdot) =  A(\cdot,\by,\bx,\cdot) \text{ in }(0,T)\times\R, \text{ for }\bx,\by\in\O.
\end{align}
\item The reaction function $f \in L^\infty((0,T)\times \O) \times W_{loc}^{1,\infty}(\R)$ satisfies 
\begin{align}
& {\redtext QUITAR ~~~f(\cdot,\cdot,0)\geq 0 \text{ in }Q_T,} \nonumber\\
& \abs{f(\cdot,\cdot,s)}\leq C_f (1+\abs{s}) \text{ in }Q_T, \text{ for }s\in\R, \text{ and constant } C_f>0. \label{cond.f2}
\end{align}

\item The initial datum $u_0\in L^1(\O)$ is non-negative. 
\end{enumerate}  

Before stating our results, let us introduce the notion of solution we employ for problem \fer{eq.eq}-\fer{eq.id}.
We interpret both equations in the a.e. pointwise sense. 
\begin{definition}
A solution of problem \fer{eq.eq}-\fer{eq.id} is a function $u\in W^{1,1}(0,T;L^1(\O))$ such that 
\begin{align*}
&  \p_t u(t,\bx)  =  \int_\O J(\bx-\by) A (t,\bx,\by,u(t,\by)-u(t,\bx)) d\by +f(t,\bx,u(t,\bx)), 
\end{align*}
for a.e. $(t,\bx)\in Q_T$, and with $u(0,\cdot)=  u_0$ a.e. in $\O$. 
\end{definition}

Our first result states the well-posedness of problem \fer{eq.eq}-\fer{eq.id} in the case of regular data.
The spatial regularity  required in \fer{ass:add0}  is due to the integral relationship that $J$ and $A$ play on equation \fer{eq.eq}. Thus, it can be somehow varied by, for instance, weakening the regularity demanded for $A$ at the cost of requiring more regularity for $J$. We comment on this variations at the end of this section. In particular, Theorem~\ref{th.existence-Lipschitz} ensures the existence of solution of the gradient dependent nonlocal $p$-Laplacian evolution problem \fer{eq:varplap}, or to the gradient descent for the functionals corresponding to bilateral filters of the type \fer{def.bil} and  \fer{def.yar}.

\begin{theorem}
 \label{th.existence-Lipschitz}
Assume (H) and, additionally, 
\begin{align}
& J \in BV(B),&&{\redtext A \in L^\infty(0,T;W^{1,\infty}(\O\times\O)\times W_{loc}^{1,\infty}(\R))}, & \label{ass:add0}\\
& {\redtext f(t,\cdot,s)\in BV(\O)},&&  u_0\in  L^\infty(\O)\cap BV(\O) ,& \label{ass:add1}
\end{align}
for $t\in(0,T)$ and $s \in \R$.
Then there exists a unique solution of problem \fer{eq.eq}-\fer{eq.id}, 
\begin{align*}
u\in W^{1,\infty}(0,T;L^\infty(\O))\cap C([0,T];L^\infty(\O)\cap BV(\O)),
\end{align*}
such that, for some constant $C_1>0$ depending only on $\nor{u_0}_{L^\infty(\O)}$ and $C_f$,   
\begin{align*}
 \nor{u}_{L^\infty(Q_T)} \leq C_1.
\end{align*}
{\redtext In addition, if $f(\cdot,\cdot,0)\geq 0$  in $Q_T$ then the solution is non-negative.} Finally, suppose that 
$J\in L^\infty(B)$ and let $u_1,u_2$ be the solutions of problem \fer{eq.eq}-\fer{eq.id} corresponding to the initial data $u_{01},u_{02}$ then, for a.e. $t \in(0,T)$,
\begin{equation}
\label{stability}
 \nor{u_1(t,\cdot)-u_2(t,\cdot)}_{L^\infty(\O)}\leq C_2\nor{u_{01}-u_{02}}_{L^\infty(\O)},
\end{equation}
for some constant $C_2>0$, depending only upon  $C_f,  L_A, L_f$ and $\nor{J}_{L^\infty(B)}$, where 
$L_A$ and $L_f$ are Lipschitz continuity constants for $A$ and $f$.

\end{theorem}

Our second result establishes the well-posedness of problem \fer{eq.eq}-\fer{eq.id} for non-regular data, at the cost of assuming a monotonic behavior on the range kernel, $A$. This result ensures the existence of solution
for the nonlocal $p(t,\bx,\by)-$Laplacian evolution problem, with $p\in(1,\infty)$. 

\begin{theorem}
 \label{th.approximation}
Assume (H) and suppose that $A(t,\bx,\by,\cdot)$ is non-decreasing in $\R$, for $t\in(0,T)$ and for $\bx,\by\in\O$.  Let $u_0\in L^q(\O)$ for some $q\in[1,\infty]$, and assume that
\begin{align*}
J\in L^{\frac{q}{q-\alpha}}(B) \text{ if }q\in[1,\infty) \qtext{or}\quad J\in L^1(B) \text{ if }q=\infty.
\end{align*}
Then, there exists a unique solution $u\in W^{1,1}(0,T;L^1(\O))\cap L^q(Q_T)$ of problem \fer{eq.eq}-\fer{eq.id}.
{\redtext In addition, if $f(\cdot,\cdot,0)\geq 0$  in $Q_T$ then the solution is non-negative.}  Moreover, if $u_1,u_2$ are solutions corresponding to the initial data $u_{01},u_{02}\in L^q(\O)$ then, for a.e. $t\in(0,T)$,
\begin{equation}
 \label{stability2}
 \nor{u_1(t,\cdot)-u_2(t,\cdot)}_{L^q(\O)}\leq C\nor{u_{01}-u_{02}}_{L^q(\O)},
\end{equation}
for some constant $C>0$. Besides, $C=1$ if $f$ is non-increasing.
\end{theorem}

Some extensions and variations are possible for the hypothesis assumed in Theorems~\ref{th.existence-Lipschitz} and \ref{th.approximation}.  We have omitted them to keep a reasonable clarity in the exposition. 
However, the proofs may be easily modified to incorporate them. We list here some possibilities.

For both theorems, we may directly consider a symmetric space kernel $J:\R^d\times\R^d\to\R$, that is, satisfying $J(\bx,\by)=J(\by,\bx)$. This is in connection to the symmetry of $A$ assumed on \fer{ass:asym}, which is essential for the integral term of equation \fer{eq.eq} to be a dissipative operator, and, thus, a nonlocal diffusion, see Lemma~\ref{lemma:intparts}.
In addition, the integral condition \fer{cond.j1} is not necessary, and has been introduced just for avoiding recurrent unimportant constants arising in the estimations.

The sub-linearity assumption for $f$ given in \fer{cond.f2} is imposed to obtain global in time existence of solutions. It can be dropped for a local existence result. Besides, in Theorem~\ref{th.approximation}, 
we may replace the assumption on local Lipschitz continuity on $f$ by requiring H\"older continuity and decreasing monotony.

For the existence of solutions in both theorems, a nonlocal reaction term of the form 
\begin{align*}
\int_\O J_r(\bx-\by) F(u(t,\by)) d\by,
\end{align*}
may be included in the equation \fer{eq.eq}, with $J_r$ and $F$ satisfying similar properties than those imposed on $J$ and $f$. 

Finally, in Theorem~\ref{th.existence-Lipschitz}, we may replace the Lipschitz continuity assumed on the space variables of $A$ by $A(t,\cdot,\cdot,s) \in W^{1,\beta}(\O\times\O)$, with $\beta=d/(d-1)$ if $\O\subset\R^d$ with $d>1$, or $\beta=\infty$ if $d=1$. This is related to the optimal embedding of the space $BV$ into $L^p$ spaces. Moreover, if we further assume in that theorem that $J\in L^\infty(B)$, then the result holds for $A(t,\cdot,\cdot,s) \in BV(\O\times\O)$.

\begin{remark}
\label{remark1} The nonlocal ${p(u(t,\by)-u(t,\bx))}$-Laplacian.

Let us consider the range kernel given by 
\begin{align}
\label{pgradulap}
A(s) = \abs{s}^{p(\abs{s})-2} s,
\end{align}
for some Lipschitz continuous non-increasing $p:[0,\infty)\to\R$ satisfying, at least, one of the following conditions:
\begin{enumerate}
\item $p(0)> 2$, or 
\item $p(0)=2$ and $p'(0)<0$.
\end{enumerate}
Then, Theorem~\ref{th.existence-Lipschitz} provides the existence of a unique solution of problem \fer{eq.eq}-\fer{eq.id} for $A$ given by \fer{pgradulap}. Indeed, we have for $\sigma = \abs{s}$
\begin{align*}
A'(s) = \sigma^{p(\sigma)-2} (p(\sigma)+\sigma \log (\sigma) p'(\sigma)-1),
\end{align*}
which is bounded in bounded intervals, and thus satisfies the conditions of the theorem. Observe that the function $p$ may take values smaller than one outside $s=0$, including in this way the hyper-Laplacian diffusion \cite{Krishnan2009}.
\end{remark}

\section{Proofs}

The following integration formula, resembling integration by parts in differential calculus, is a consequence of the symmetry properties of the kernels $J$ and $A$. It also implies that  the integral term of equation \fer{eq.eq} is dissipative. The proof is straightforward, so we omit it.
\begin{lemma}
 \label{lemma:intparts}
 Assume (H) and let $u,\vfi \in L^\infty(Q_T)$, $\rho\in L^\infty(0,T)$, with $\rho \geq 0$ in $(0,T)$. Then, for a.e. $t\in(0,T)$, 
\begin{align}
\label{prop.parts}
\int_\O   \int_\O  & J(\bx-\by) A\big(t,\bx,\by,\rho(t)(u(t,\by)-u(t,\bx))\big) \vfi(t,\bx) d\by d\bx \\
 = &  - \frac{1}{2} \int_\O \int_\O J(\bx-\by) A\big(t,\bx,\by,\rho(t)(u(t,\by)-u(t,\bx))\big) \nonumber \\
 & \qquad\qquad  \times (\vfi(t,\by)-\vfi(t,\bx)) d\by d\bx. 
 \nonumber
\end{align}
In particular, for $\vfi=\phi(u)$, with $\phi\in L^\infty(\R)$ non-decreasing, we have, for a.e. $t\in(0,T)$, 
\begin{align}
\label{prop.inc}
\int_\O  & \int_\O  J(\bx-\by) A\big(t,\bx,\by,\rho(t)(u(t,\by)-u(t,\bx))\big) \phi(u(t,\bx)) d\by d\bx \leq 0.
\end{align}
\end{lemma}

\medskip

\no\textbf{Proof of Theorem~\ref{th.existence-Lipschitz}.}
We divide the proof in two steps. In the first step, we prove the result for data with more regularity than the assumed in the theorem. 
We discretize in time, find estimates for the time-independent sequence of problems, and pass to the limit in the discretization parameter with the use of compactness arguments. 
In the second step, we pass to the limit with respect to the regularized data using a similar compactness technique.

\medskip

\noindent\textbf{Step 1. Regularized problem.} 

 We first prove the existence of solutions of problem \fer{eq.eq}-\fer{eq.id} for the case in which, in addition to (H),   \fer{ass:add0} and \fer{ass:add1}, we have  
 \begin{align}
 \label{ass:add2}
J\in W^{1,1}(B), \quad  u_0\in W^{1,\infty}(\O),\quad {\redtext f(t,\cdot,s) \in W^{1,\infty}(\O)},
 \end{align}
{\redtext for $t\in(0,T)$, and $s\in\R$. We also assume the following sub-linearity condition on the range kernel,
 \begin{align}
  \abs{A(\cdot,\cdot,\cdot,s)}\leq C_A \abs{s}, \qtext{in } (0,T)\times\overline{\O}\times\overline{\O}, \qtext{for }s\in\R. \label{cond.a1}
 \end{align}
 for some constant $C_A>0$. This condition is instrumental to this step. It will be used to obtain preliminar $L^\infty$ estimates for the solutions of some time semi-discrete approximated problems. After passing to the limit in the time discretization, the $L^\infty$ estimate of the emerging solution is shown to be independent of $C_A$.}

Consider the following auxiliary problem, obtained using the change of unknown $u=\e^{\mu t} w$ in \fer{eq.eq}, for some positive constant $\mu$ to be fixed:
\begin{align}
 \p_t w(t,\bx)  =& ~ \e^{-\mu t} \int_\O J(\bx-\by) A \big(t,\bx,\by,\e^{\mu t}(w(t,\by)-w(t,\bx))\big) d\by \nonumber \\ 
 & + \e^{-\mu t} f(t,\bx,\e^{\mu t}w(t,\bx)) - \mu w(t,\bx), \label{eq.aux} \\
 w(0,\bx) =& ~u_0(\bx), \label{id.aux}
 \end{align}
 for $(t,\bx)\in (0,T)\times\O$. 

\smallskip 

\no\emph{Time discretization.} Let $N\in \N$, $\tau=T/N$ and $t_j=j\tau$, for $j=0,\ldots,N$. 
Assume that $w_j\in W^{1,\infty}(\O)$ is given and consider the following time discretization of \fer{eq.aux}:
 \begin{align}
 w_{j+1}(\bx)  = & w_{j}(\bx)+\tau \e^{-\mu t_j} \int_\O J(\bx-\by) A \big(t_j,\bx,\by,\text{e}^{\mu t_j}(w_{j}(\by)-w_{j}(\bx))\big) d\by   \nonumber \\
 & +\tau  \e^{-\mu t_j} f(t_j,\bx,\e^{\mu t_j}w_j(\bx))  - \tau \mu w_{j+1}(\bx). \label{eq.discrete}
 \end{align}

\smallskip 

\no\emph{Uniform estimates with respect to $\tau$. } 
Let us show that $w_{j+1}$ is uniformly bounded in $W^{1,\infty}(\O)$. On one hand, using in \fer{eq.discrete} the growth conditions \fer{cond.a1} and \fer{cond.f2} on $A$ and $f$, together with the normalization property \fer{cond.j1} on $J$, we obtain 
\begin{align*}
 \nor{w_{j+1}}_{L^\infty} & \leq \frac{1}{1+\tau \mu}\Big(\nor{w_j}_{L^\infty} +\tau (2 C_A + C_f)  \nor{w_j}_{L^\infty} +\tau C_f\Big) \\
& = \frac{1+\tau (2 C_A+C_f)}{1+\tau \mu}\nor{w_j}_{L^\infty} +\frac{\tau C_f}{1+\tau \mu}. 
\end{align*}
Taking $\mu > 2C_A + C_f$, this differences inequality yields the uniform estimate 
\begin{equation}
 \label{est.w}
 \nor{w_{j+1}}_{L^\infty} \leq M_0,
\end{equation}
with $M_0$ depending on the regularity assumed on this step only through $C_A$.

On the other hand, taking into account the assumptions \fer{ass:add1} and \fer{ass:add2},  and the estimate \fer{est.w}, we deduce from \fer{eq.discrete} that  $w_{j+1}\in W^{1,\infty}(\O)$. This regularity allows to 
differentiate  in \fer{eq.discrete} with respect to the $k-$th component of $\bx$, denoted by $x_k$, to obtain for a.e. $\bx\in\O$,
\begin{align*}
 (1+\tau \mu) \frac{\p w_{j+1}}{\p x_k}(\bx) = F_1(\bx) \frac{\p w_{j}}{\p x_k}(\bx) + \tau \e^{-\mu t_j} F_2(\bx) , 
 \end{align*}
 with 
\begin{align*}
F_1(\bx)  = & 1 - \tau \int_\O J(\bx-\by) \frac{\p A}{\p s} \big(t_j,\bx,\by,\text{e}^{\mu t_j}(w_{j}(\by)-w_{j}(\bx)) \big) d\by  \\
&+ \tau  \frac{\p f}{\p s}(t_j,\bx, \e^{\mu t_j}w_j(\bx)) ,  \\
 F_2(\bx)  = &  \int_\O \frac{\p J}{\p x_k}(\bx-\by) A \big(t_j,\bx,\by,\text{e}^{\mu t_j}(w_{j}(\by)-w_{j}(\bx)) \big) d\by  \\
 + & \int_\O J(\bx-\by) \frac{\p A}{\p x_k} \big(t_j,\bx,\by,\text{e}^{\mu t_j}(w_{j}(\by)-w_{j}(\bx)) \big) d\by  \\
 & + \frac{\p f}{\p x_k}(t_j,\bx, \e^{\mu t_j}w_j(\bx)) .
  \end{align*}
We deduce
\begin{align*}
  \nor{\nabla w_{j+1}}_{L^\infty} \leq & \frac{1}{1+ \tau \mu } \Big(  \big(1+ \tau (L_A+L_f)  \big)  \nor{\nabla w_j}_{L^\infty} + \tau 2L_A \nor{\grad J}_{L^1}  \nor{w_j}_{L^\infty}\\
  & +\tau (L_A+L_f)\Big), 
 \end{align*}
 where $L_A$ is the Lipschitz constant of $A(t_j,\cdot,\cdot,\cdot)$ in $\O\times\O\times [-2\e^{\mu T}M_0,2\e^{\mu T}M_0]$, and $L_f$ is the Lipschitz constant of $f(t_j,\cdot,\cdot)$ in $\O\times[-\e^{\mu T}M_0,\e^{\mu T}M_0]$. Choosing 
\begin{align*}
\mu > \max\{2C_A+C_f,  2M_0L_A \nor{\grad J}_{L^1} +L_A+L_f \},
\end{align*}
and solving this differences inequality, we obtain the uniform estimate 
\begin{align}
\label{est.gradw}
\nor{\nabla w_{j+1}}_{L^\infty}\leq M_1,
\end{align}
with $M_1$ depending on the regularization introduced in this step only through  $\nor{J}_{W^{1,1}}$, $L_A$,  $C_A$, and $\nor{u_0}_{W^{1,\infty}}$. 

\bigskip 
\no\emph{Time interpolators and passing to the limit $\tau\to0$. }

We define, for $(t,\bx)\in (t_j,t_{j+1}]\times \O$, the piecewise constant and
piecewise linear interpolators of $w_j$ given by
\begin{align*}
 w^{(\tau)}(t,\bx)=w_{j+1}(\bx),\quad  
 \tilde w^{(\tau)}(t,\bx)=w_{j+1}(\bx)+\frac{t_{j+1}-t}{\tau}(w_j(\bx)-w_{j+1}(\bx)),
\end{align*}
and the piecewise constant  approximations 
\begin{align*}
  \text{e}_\tau^{\pm\mu t}=\text{e}^{\pm\mu t_j}, \quad A_\tau (t,\cdot,\cdot,\cdot) = A(t_j,\cdot,\cdot,\cdot), \quad f_\tau (t,\cdot,\cdot) = f(t_j,\cdot,\cdot),
\end{align*}
which converge in $L^p(0,T)$, for any $p\in[1,\infty)$, and pointwise for a.e. $t\in(0,T)$ to the exponential function, to $A(t,\cdot, \cdot,\cdot)$, and to $f(t,\cdot,\cdot)$, respectively.

We also introduce the shift operator $\sigma_\tau w^{(\tau)}(t,\cdot) = w_{j}$. With this notation,  equation \fer{eq.discrete} may be rewritten as, for $(t,\bx)\in Q_T$,  
\begin{align}
 \p_t \tilde w^{(\tau)}(t,\bx)  =& \e_\tau^{-\mu t} \int_\O J(\bx-\by) A_\tau \big(t,\bx,\by,\e_\tau^{\mu t}(\sigma_\tau w^{(\tau)}(t,\by)-\sigma_\tau w^{(\tau)}(t,\bx))\big) d\by \nonumber \\ 
 & + \e_\tau^{-\mu t} f_\tau(t,\bx,e_\tau^{\mu t}\sigma_\tau w^{(\tau)}(t,\bx))  - \mu w^{(\tau)}(t,\bx). \label{eq.interp}
 \end{align}

 Using the uniform $L^\infty$ estimates \fer{est.w} and \fer{est.gradw} of $w_{j+1}$ and $\nabla w_{j+1}$, we deduce the corresponding uniform estimates for $\nor{ \nabla w^{(\tau)}}_{L^{\infty}(Q_{T})}$, 
$\nor{ \nabla \tilde w^{(\tau)}}_{L^{\infty}(Q_{T})}$ and
$\nor{ \p_t \tilde w^{(\tau)}}_{L^{\infty}(Q_{T})}$, implying the existence of 
$w\in L^\infty(0,{T}; W^{1,\infty}(\O))$ and 
$\tilde w\in W^{1,\infty}(Q_{T})$ such that, at least for subsequences (not relabeled)
\begin{align}
 & w^{(\tau)} \to w \qtext{weakly* in }L^\infty(0,T; W^{1,\infty}(\O)), \nonumber \\
 & \tilde w^{(\tau)} \to \tilde w \qtext{weakly* in }   W^{1,\infty}(Q_{T}), \label{conv.1}
\end{align}
as $\tau \to 0$.
In particular, by compactness 
\begin{equation*}
 \tilde w^{(\tau)} \to \tilde w \qtext{uniformly in } C([0,T]\times \bar \O). 
\end{equation*}
Since, for $t\in (t_j,t_{j+1}]$,
\begin{align*}
 \abs{ w^{(\tau)}(t,\bx) -\tilde w^{(\tau)}(t,\bx) }=& \left|\frac{(j+1)\tau -t}{\tau}(w_j(\bx)-w_{j+1}(\bx))\right| \\
 & \leq \tau \nor{\p_t \tilde w^{(\tau)}}_{L^{\infty}(Q_{T})},
\end{align*}
we deduce both $w=\tilde w$ and, up to a subsequence, 
\begin{equation}
 \label{conv.2}
 w^{(\tau)} \to  w \qtext{strongly in } L^\infty(Q_T) \text{ and a.e. in } Q_T. 
\end{equation}
With the properties of convergence \fer{conv.1} and \fer{conv.2} the passing to the limit $\tau\to0$ in \fer{eq.interp} is justified, finding that $w\in W^{1,\infty}(Q_T)$ is a solution of  \fer{eq.aux}-\fer{id.aux}, and therefore, $u=w\text{e}^{\mu t} \in W^{1,\infty}(Q_T)$  is a solution of \fer{eq.eq}-\fer{eq.id}.

 \medskip
 
\noindent\emph{A posterior estimates. } Using properties \fer{prop.parts} and \fer{prop.inc}, we show uniform estimates with respect to the local Lipschitz continuity and the growth condition on $A$ (constants $L_A$ and $C_A$), and to the  Lipschitz continuity of $f(t,\bx,\cdot)$ (constant $L_f$). 

(i) $L^\infty$ bound. {\redtext We show that the positive part of $u$, denoted by $u_+$, is bounded in $L^\infty(Q_T)$. The result for the negative part is achieved similarly.}
Consider again the change of unknown $u= e^{\mu t}w$ leading to equation \fer{eq.aux}, and multiply this equation by $\phi(w)\in W^{1,\infty}(Q_T)$, for the non-decreasing function $\phi(s)=\max\{0,s-K\}$, with $K>0$ to be fixed. Using Lemma~\ref{lemma:intparts}, we get  
\begin{align*}
  \p_t \int_\O \Phi(w(t,\bx))d\bx \leq  \int_\O \big(\e^{-\mu t} f(t,\bx,\e^{\mu t}w(t,\bx)) - \mu w(t,\bx)\big)\phi(w(t,\bx))d\bx, 
 \end{align*}
where $\Phi(s)=\phi(s)^2/2$. The growth condition on $f(t,\bx,\cdot)$  implies
\begin{align*}
  \p_t \int_\O \Phi(w(t,\bx))d\bx \leq & (C_f-\mu) \int_\O \phi(w(t,\bx)) \abs{w(t,\bx)- K}    d\bx \\
&   +  \int_\O \big(C_f (\e^{-\mu t} +K) -\mu K \big) \phi(w(t,\bx))d\bx .
 \end{align*}
Taking $\mu = C_f(1+K)/K$, and noticing that $\phi$ is non-negative, we deduce 
\begin{align*}
  \p_t \int_\O \Phi(w(t,\bx))d\bx \leq (C_f-\mu) \int_\O \abs{\phi(w(t,\bx))}^2   d\bx. 
 \end{align*}
Fixing $K> \nor{u_0}_{L^\infty}$, Gronwall's lemma yields, for the original unknown,
\begin{align}
 \label{est.linf.unif}
 \nor{u_+}_{L^\infty(Q_T)} \leq e^{\mu T}\nor{u_0}_{L^\infty(\O)},
\end{align}
with $\mu$ depending only on $C_f$ and $\nor{u_0}_{L^\infty}$.

(ii) Non-negativity. {\redtext Assume $f(\cdot,\cdot,0)\geq 0$ in $Q_T$.} We multiply \fer{eq.eq}  by $\phi(u)\in W^{1,\infty}(Q_T)$, for the non-decreasing function $\phi(s)=\min\{0,s\}$, and use Lemma~\ref{lemma:intparts} to get, for $t\in(0,T)$, 
\begin{align*}
 \p_t \int_\O \Phi(u(t,\bx))d\bx \leq \int_\O f(t,\bx,u(t,\bx))\phi(u(t,\bx))d\bx,
\end{align*}
where $\Phi(s)=\int_0^s\phi(\sigma)d\sigma = \phi(s)^2/2$. Since $f(\cdot,\cdot,0)\geq 0$ in $Q_T$, using that $\phi \leq 0$, $u\in L^\infty(Q_T)$,  and the Lipschitz continuity of $f(t,\bx,\cdot)$ we get 
\begin{align*}
 \p_t \int_\O \Phi(u(t,\bx))d\bx \leq  L_f \int_\O \abs{u(t,\bx)}\abs{\phi(u(t,\bx))}d\bx = L_f\int_\O \abs{\phi(u(t,\bx))}^2d\bx.
\end{align*}
 Therefore, Gronwall's lemma implies 
\begin{align*}
 \int_\O \abs{\phi(u(t,\bx))}^2d\bx \leq \e^{2L_f t} \int_\O \abs{\phi(u_0(\bx))}^2d\bx = 0,
\end{align*}
so that $u\geq0$ a.e. in $Q_T$.

\bigskip

\noindent\textbf{Step 2. Passing to the limit in the regularization} 

We consider sequences of smooth approximating functions $J_\eps$ and $u_{0\eps}$ satisfying 
assumptions (H),   \fer{ass:add0},  \fer{ass:add1} and \fer{ass:add2}
such that, as $\eps\to0$,
\begin{align}
 & J_{\eps} \to J \qtext{strongly in }L^1(B),  \text{ with }\nor{\nabla J_{\eps}}_{L^1(B)}\to \text{TV}(J),&&& \label{propeps0}\\
 & u_{0\eps} \to u_0 \qtext{strongly in }L^q(\O), ~\nor{u_{0\eps}}_{L^\infty(\O)}\leq K, &&&  \label{propepsinf} \\
& {\redtext f_{\eps}(t,\cdot,s) \to f(t,\cdot,s) \qtext{strongly in }L^q(\O), ~\nor{f(t,\cdot,s)}_{L^\infty(\O)}\leq K, }&&&  \label{propepsf}
\end{align}
{\redtext for $t\in(0,T)$ and $s\in\R$, and} for any $q\in[1,\infty)$, where $K>0$ is independent of $\eps$, and such that 
\begin{align}
 & \nor{\nabla u_{0\eps}}_{L^1(\O)}\to \text{TV}(u_0), \label{propeps1} \\
 & \nor{\nabla f_{\eps}(t,\cdot,s)}_{L^1(\O)}\to \text{TV}(f_{\eps}(t,\cdot,s)), \label{propepsf2} 
\end{align}
where TV denotes the total variation with respect to the $\bx$ variable. 

Sequences with properties \fer{propeps0}-{\redtext \fer{propepsf2}} do exist thanks to the regularity $J\in BV(B)$,  $u_0 \in L^\infty(\O)\cap BV(\O)$, {\redtext and $f_{\eps}(t,\cdot,s) \in BV(\O)$}, see \cite{Ambrosio2000}.
In addition, we consider a sequence $A_\eps$ with the same regularity as $A$, see \fer{ass:add0},  and 
satisfying the sub-linearity condition \fer{cond.a1}, which is possible because $A$ is locally Lipschitz continuous
in the fourth variable.  

Due to the above convergences, we have 
\begin{align}
& \grad J_{\eps}  \qtext{is uniformly bounded in }L^1(B) , \label{propeps2}\\
& \grad u_{0\eps}  \qtext{is uniformly bounded in }L^1(\O) . \label{propeps1b} \\
& {\redtext \grad f_{\eps}(t,\cdot,s)  \qtext{is uniformly bounded in }L^1(\O) . }\label{propepsf2b} 
\end{align}
Because of \fer{est.linf.unif} and \fer{propepsinf}, the corresponding solution $u_\eps\in W^{1,\infty}(Q_T)$ of problem \fer{eq.eq}-\fer{eq.id}, ensured by Step 1 of this proof, is uniformly bounded in $L^\infty(Q_T)$, 
\begin{align}
\label{est:ueinf}
\nor{u_\eps}_{L^\infty(Q_T)} \leq C,
\end{align}
for some $C>0$ independent of the growth condition constant of $A_\eps$, $C_{A_\eps}$. In particular, this means that we may replace $A_\eps$ by $A$, since the sub-linearity condition 
\fer{cond.a1} is trivially satisfied by the Lipschitz continuous function $A(\cdot,\cdot,\cdot,s)$ for $\abs{s}\leq C$. 
Therefore,  $u_\eps \in W^{1,\infty}(Q_T)$ satisfies
\begin{align}
  \p_t u_\eps(t,\bx)  = &  \int_\O J_\eps(\bx-\by) A (t,\bx,\by,u_\eps(t,\by)-u_\eps(t,\bx)) d\by \nonumber \\
&   + {\redtext f_\eps \big(t,\bx,u_\eps(t,\bx)\big)} , \label{eq.eqeps}\\
    u_\eps(0,\bx)  = & ~ u_{0\eps}(\bx),  \label{eq.ideps}
\end{align}
 for $(t,\bx)\in Q_T$.
 From \fer{est:ueinf}, the uniform boundedness of $J_\eps$ in $L^1(B)$   and  the regularity of $A$ and $f$,  we obtain from \fer{eq.eqeps} that
\begin{equation}
 \label{bound.1}
\partial_t u_\eps \qtext{is uniformly bounded in }L^\infty(Q_T). 
\end{equation}
Being $u_{0\eps}, J_\eps$ smooth functions, we may deduce an $L^\infty$ bound for $\grad u_\eps$ as in the Step~1, not necessarily uniform in $\eps$, but  allowing to differentiate equation \fer{eq.eqeps} with respect to $x_k$. After integration in $(0,t)$,  we obtain
  \begin{align}
 \label{dsv}
  \frac{\p u_\eps}{\p x_k} (t,\bx)=  G_\eps(t,\bx) \Big(  \frac{\p u_{0\eps}}{\p x_k}(\bx) +  
 \int_0^t  \eta^\eps(s,\bx) (G_\eps(s,\bx))^{-1}  ds \Big),
 \end{align}
 with 
 \begin{align*}
 G_\eps(t,\bx)  = &  \exp\Big( \int_0^t \Big[{\redtext \frac{\p f_\eps}{\p s}}(t,\bx,u_\eps(\tau,\bx)) \nonumber \\
 & \qquad - \int_\O J_\eps(\bx-\by) \frac{\p A}{\p s} \big(t,\bx,\by, u_\eps(\tau,\by)-u_\eps(\tau,\bx)\big)  d\by\Big] d\tau \Big), \nonumber\\
  \eta^\eps_{k}(t,\bx) =& \int_{\O} \frac{\partial J_\eps}{\partial x_k}(\bx-\by) A\big(t,\bx,\by, u_\eps(t,\by)-u_\eps(t,\bx)\big) d\by \nonumber \\
  & +\int_{\O} J_\eps(\bx-\by) \frac{\partial A}{\partial x_k}\big(t,\bx,\by, u_\eps(t,\by)-u_\eps(t,\bx)\big) d\by \nonumber \\
  & + {\redtext \frac{\p f_\eps}{\p x_k}}(t,\bx,u_\eps(\tau,\bx)). 
 \end{align*}
 Using the regularity of $A$ and properties \fer{propeps2}, \fer{propepsf2b}, and \fer{est:ueinf}, we deduce that $G_\eps$, $1/G_\eps$  are uniformly bounded in $L^\infty(Q_T)$, and that $\eta^\eps$ is uniformly bounded in $L^\infty(0,T;L^1(\O))$. Then, from \fer{dsv} and the uniform bound \fer{propeps1b},  we obtain  that 
\begin{equation}
 \label{bound.2}
\nabla u_\eps \qtext{is uniformly bounded in }L^\infty(0,T;L^1(\O)). 
\end{equation}
Bounds \fer{bound.1} and \fer{bound.2} allow to deduce, using the compactness result \cite[Cor. 4, p. 85]{Simon1986},  the existence of $u\in C([0,T];L^\infty(\O)\cap BV(\O))$ such that $u_\eps\to u$ strongly in  $L^q(Q_T)$, for all $q<\infty$, and a.e. in $Q_T$. 
The uniform bound \fer{bound.1} also implies that, up to a subsequence (not relabeled), we have 
$\p_t u_\eps \to \p_t u$ weakly* in $L^\infty(Q_T)$. 

These convergences allow to pass to the limit $\eps\to 0$ in \fer{eq.eqeps}-\fer{eq.ideps} 
(with $u$ replaced by $u_\eps$) and identify the limit 
\begin{align*}
u \in W^{1,\infty}(0,T;L^\infty(\O))\cap C([0,T];L^\infty(\O)\cap BV(\O)), 
\end{align*}
 as a solution of \fer{eq.eq}-\fer{eq.id}. 
Observe that, being $u_\eps$ a sequence of non-negative functions, we deduce $u\geq 0$ in $Q_T$.
  
\smallskip  
  
\noindent\emph{Uniqueness and stability. }
Let $u_{01},~u_{02}\in L^\infty(\O)\cap BV(\O)$ and $u_1,~u_2$ be the corresponding solutions to problem
\fer{eq.eq}-\fer{eq.id}. Set $u=u_1-u_2$ and   $u_0=u_{10}-u_{20}$. Then $u \in W^{1,\infty}(0,T;L^\infty(\O))\cap C([0,T];L^\infty(\O)\cap BV(\O))$ satisfies $u(0,\cdot)= u_0$ in $\O$, and
\begin{align}
 \p_t u(t,\bx)  = 
 \int_\O J(\bx-\by) & \Big( A(t,\bx,\by,u_1(t,\by)-u_1(t,\bx)) \nonumber \\
 & - A(t,\bx,\by,u_2(t,\by)-u_2(t,\bx))\Big)d\by \nonumber \\
& + f(t,\bx,u_1(t,\bx))- f(t,\bx,u_2(t,\bx)), \label{eq.stabi}
 \end{align}
 for $(t,\bx)\in Q_T$. 
  Multiplying this equation by $u$, integrating in $\O$ and using the  Lipschitz continuity of $A(t,\bx,\by,\cdot)$ and $f(t,\bx,\cdot)$, we deduce
\begin{align*}
  \frac{1}{2}\p_t \int_\O \abs{u(t,\bx)}^2 d\bx \leq & L_A \int_\O \int_\O J(\bx-\by)  \abs{u(t,\by)-u(t,\bx)} \abs{u(t,\bx)}d\by d\bx  \\
& + L_f\int_\O \abs{u(t,\bx)}^2d\bx.
\end{align*}
Using the inequality $\abs{s}\abs{t-s}\leq 2 (t^2+s^2)$ and the summability property of $J$ in the first term of the right hand side, we obtain
\begin{align*}
  \frac{1}{2}\p_t \int_\O \abs{u(t,\bx)}^2 d\bx \leq & (2L_A +L_f)\int_\O \abs{u(t,\bx)}^2d\bx.
\end{align*}
Therefore, if $u_{10}=u_{20}$, Gronwall's inequality implies $u_1=u_2$.

For the stability result, we assumed $J\in L^\infty(B)$.  
 Multiplying \fer{eq.stabi} by $\phi(u)$, with $\phi(s)=\abs{s}^{q-1}s$, for $q\geq 1$, integrating in $\O$ and using the  Lipschitz continuity of $A(t,\bx,\by,\cdot)$ and $f(t,\bx,\cdot)$, and the boundedness of $J$, we deduce
\begin{align*}
 \frac{1}{q+1}  \p_t \int_\O \abs{u(t,\bx)}^{q+1} d\bx \leq & L_A \nor{J}_{L^\infty}\int_\O \int_\O  \abs{u(t,\by)-u(t,\bx)} \abs{u(t,\bx)}^{q}d\by d\bx  \\
& + L_f\int_\O \abs{u(t,\bx)}^{q+1}d\bx.
\end{align*}
By Young's inequality we find, for fixed $t\in(0,T)$,  
\begin{align*}
\int_\O \int_\O  \abs{u(t,\by)} \abs{u(t,\bx)}^{q}d\by d\bx = \nor{u}_{L^1}\nor{u}_{L^q}^q \leq \abs{\O}^{\frac{2q-1}{q}}
\nor{u}_{L^{q+1}}^{q+1}.
\end{align*}
Therefore, we deduce  
\begin{align*}
   \p_t \int_\O \abs{u(t,\bx)}^{q+1} d\bx \leq & (q+1)\Big(L_f + L_A \nor{J}_{L^\infty} ( \abs{\O}^{\frac{2q-1}{q}} +\abs{\O})\Big)  \int_\O \abs{u(t,\bx)}^{q+1}d\bx,
\end{align*}
and Gronwall's inequality implies, for a.e. $t\in(0,T)$, 
\begin{align*}
\nor{u(t,\cdot)}_{L^{q+1}(\O)} \leq \exp\Big(t \big(L_f + L_A \nor{J}_{L^\infty} ( \abs{\O}^{\frac{2q-1}{q}} +\abs{\O})\big)\Big)\nor{u_0}_{L^{q+1}(\O)}.
\end{align*}
The result follows letting $q\to\infty$.
$\Box$

\medskip

In the proof of Theorem~\ref{th.approximation} we use the following lemmas: an approximation result, and a consequence of the monotonicity of $A$. We prove them at the end of this section.
\begin{lemma}
 \label{lemma:mona}
Assume (H) and suppose that $A(t,\bx,\by,\cdot)$ is non-decreasing in $\R$, for $t\in(0,T)$ and for $\bx,\by\in\O$.  Then, there exists a sequence 
$A_n \in L^\infty(0,T)\times W^{1,\infty}(\O\times\O)\times W_{loc}^{1,\infty}(\R)$, for $n\in\N$,  such that 
\begin{align}
& A_n(\cdot,\cdot,\cdot,0)= 0 \text{ in }(0,T)\times\overline{\O}\times\overline{\O},  \text{ for all }n\in\N ,\label{p.a0}\\
& \frac{dA_n}{ds} (\cdot,\cdot,\cdot,s)\geq 0 \text{ in }(0,T)\times\overline{\O}\times\overline{\O},  \text{ for all }s\in\R,~n\in\N \label{p.a1}\\
&  \abs{A_n(\cdot,\cdot,\cdot,s_1)  - A_n(\cdot,\cdot,\cdot,s_2) } \leq C_0 \abs{s_1-s_2}^\alpha  \text{ in }(0,T)\times\overline{\O}\times\overline{\O}, \nonumber \\
& \qquad \text{ for all } s_1,s_2\in\R,~ n\in\N, \label{p.a4}\\
 & \abs{A_{n_1}  - A_{n_2} } \leq C_0 \Big(\frac{\abs{n_1-n_2}}{n_1 n_2}\Big)^\alpha \text{ in }(0,T)\times\overline{\O}\times\overline{\O}\times\R, \nonumber \\
 & \qquad \text{ for all }n_1,n_2\in\N, \label{p.a2}\\
&  A_n \to A \text{ in } L^q((0,T)\times \O\times\O)\times L^\infty(\R)\text{ as }n\to\infty,  \label{p.a3}
\end{align}
for any $q\in[1,\infty)$, with $\nor{A_n}_{L^\infty} \leq K$, 
where $C_0, K$ are positive constants independent of $n$.
\end{lemma}

\begin{lemma}
\label{lemma:mon}
Assume (H) and suppose that $A(t,\bx,\by,\cdot)$ is non-decreasing in $\R$, for $t\in(0,T)$ and for $\bx,\by\in\O$.  
 Let $\phi:\R\to\R$ be  non-decreasing and 
$\cL$ be defined by
\begin{align*}
\cL(t) u(\bx) = -\int_\O J( \bx-\by) A(t,\bx,\by,u(\by)-u(\bx)) d\by,\qtext{for }t\in(0,T).
\end{align*}
Suppose that $\phi(u-v) (\cL(t) u-\cL(t) v)\in L^1(\O)$ for $t\in(0,T)$. Then 
\begin{align*}
  \int_\O \phi\big(u(\bx)-v(\bx)\big) (\cL(t) u(\bx)-\cL(t) v(\bx))d\bx \geq 0\qtext{for }t\in(0,T).
\end{align*}
In particular, the conditions of this lemma are satisfied if $\phi\in L^\infty(\R)$, $u,v\in L^q(\O)$, for $q\in[1,\infty]$,   and $J\in L^{\frac{q}{q-\alpha}}(\R)$  if $q\in[1,\infty)$ or $J\in L^1(B)$  if $q=\infty$.
\end{lemma}

\no\textbf{Proof of Theorem~\ref{th.approximation}.}

We consider the sequence $A_i$, for $i\in\N$,  provided by Lemma~\ref{lemma:mona}, and other sequences $J_i \in  BV(B)$, and $u_{0i}\in L^\infty(\O)\cap BV(\O)$ such that, if $q\in[1,\infty)$ 
\begin{align}
& J_{i}\to J \qtext{strongly in }L^{\frac{q}{q-\alpha}}(B), \nonumber \\
& u_{0i}\to u_0 \qtext{strongly in }L^q(\O), \label{and2}
\end{align}
or, if $q=\infty$,
\begin{align}
& J_{i}\to J \qtext{strongly in }L^1(B), \nonumber \\
& u_{0i}\to u_0 \qtext{strongly in }L^r(\O), \text{ for }r\in[1,\infty), \text{ with }\nor{u_{0i}}_{L^\infty(\O)} < C, \label{and2b}
\end{align}
for some $C>0$ independent of $i$. We  set the problem 
\begin{align}
 & \p_t u(t,\bx)  =  \int_\O J_i(\bx-\by) A_i (t,\bx,\by,u(t,\by)-u(t,\bx)) d\by +f (t,\bx,u(t,\bx)), \label{eq.eq.ap}\\
&  u(0,\bx)= u_{0i}(\bx),  \label{eq.id.ap}
\end{align}
for $(t,\bx)\in Q_T$, for which the existence of a unique solution
\begin{align*}
u_i\in W^{1,\infty}(0,T;L^\infty(\O))\cap C([0,T];L^\infty(\O)\cap BV(\O)),
\end{align*}
is  ensured by Theorem~\ref{th.existence-Lipschitz}.

We start proving the uniform boundedness of $u_i$ in $L^q(Q_T)$.  To do this, we modify the argument employed  for proving the stability result \fer{stability} of Theorem~\ref{th.existence-Lipschitz} to take into account the monotonicity of $A_i$. 

Taking $J=J_i$, $A=A_i$, $u_1=u_i$ and $u_2=0$ in \fer{eq.stabi} and using $\phi(u)$, with $\phi(s)=\abs{s}^{r-1}s$, for $r\geq 1$, as a test function in \fer{eq.stabi} we find, due to the monotonicity of $A_i$  and $\phi$, see Lemma~\ref{lemma:mon},  
and to the Lipschitz continuity of $f$,
\begin{align*}
 \frac{1}{r+1}  \p_t \int_\O \abs{u_i(t,\bx)}^{r+1} d\bx \leq L_f\int_\O \abs{u_i(t,\bx)}^{r+1}d\bx.
\end{align*}
 Then,  Gronwall's inequality implies, for any $r\in[1,\infty]$,
\begin{align*}
\nor{u_i}_{L^r(Q_T)} \leq \e^{L_fT} \nor{u_{0i}}_{L^r(\O)}.
\end{align*}
Due to \fer{and2} or \fer{and2b}, we deduce that  $\nor{u_i}_{L^q(Q_T)}$, for whatever the choice of $q\in[1,\infty]$, is uniformly bounded with respect to $i$. In particular, we have that the integral term in \fer{eq.eq.ap} is well defined for all $i \in\N$ since the H\"older continuity of $A_i$ and H\"older's inequality imply
\begin{align*}
\int_\O J_i(\bx-\by) A_i (t,\bx,\by,u_i(t,\by)-u_i(t,\bx)) d\by \leq C \nor{J_i}_{L^\frac{q}{q-\alpha}(\R^d)} \nor{u_i(t,\cdot)}_{L^q(\O)}^\alpha.
\end{align*}
This bound makes sense for $q\in[1,\infty)$. In the case $q=\infty$ we just replace $q/(q-\alpha)$ by $1$. 
In order to treat both cases jointly, we introduce the notation
\begin{align*}
\frac{q}{q-\alpha} = 
\begin{cases}
\frac{q}{q-\alpha} & \text{if }q\in[1,\infty)\\
1 & \text{if } q = \infty.
\end{cases}
\end{align*}

The main ingredient of the proof is showing that $u_i$ is a Cauchy sequence in $L^1(Q_T)$.

Let $u_m,~u_n$ be the solutions of \fer{eq.eq.ap}-\fer{eq.id.ap} corresponding to $i=m$ and $i=n$, respectively. Subtracting the corresponding equations and multiplying by $\phi_\eps(u_m-u_n)$, for $\eps>0$, where $\phi_\eps \in C(\R)$ is a non-decreasing bounded approximation of the sign function, e.g. $\phi_\eps(s)=s/\eps$ if $s\in (0,\eps)$,  $\phi(s) = 1$ if $s>\eps$, and $\phi_\eps(s) = -\phi_\eps(-s)$, if $s<0$, we get 
\begin{align}
 \p_t \int_\O & \Phi_\eps(u_m(t,\bx)  -u_n(t,\bx)) d\bx = I_1 + I_2 +I_3+I_4, \label{me}
\end{align}
being $\Phi_\eps(s)=\int_0^s\phi_\eps(\sigma)d\sigma$ an approximation of the absolute value, and 
\begin{align*}
I_1 =&  \int_\O\int_\O (J_m(\bx-\by)- J_n(\bx-\by)) \phi_\eps\big(u_m(t,\bx)-u_n(t,\bx)\big)   \\ 
&  \times   A_m\big(t,\bx,\by,u_m(t,\by)-u_m(t,\bx)\big) d\by d\bx \nonumber\\
I_2 =&  \int_\O\int_\O J_n(\bx-\by) \phi_\eps\big(u_m(t,\bx)-u_n(t,\bx)\big) 
 \nonumber\\
&  \times \Big(A_m\big(t,\bx,\by,u_m(t,\by)-u_m(t,\bx)\big) - A_m\big(t,\bx,\by,u_n(t,\by)- u_n(t,\bx)\big)\Big) d\by d\bx \nonumber\\
I_3 =&  \int_\O\int_\O J_n(\bx-\by) \phi_\eps\big(u_m(t,\bx)-u_n(t,\bx)\big) 
 \nonumber\\
&  \times  \Big(A_m\big(t,\bx,\by,u_n(t,\by)-u_n(t,\bx)\big) - A_n\big(t,\bx,\by,u_n(t,\by)- u_n(t,\bx)\big)\Big) d\by d\bx \nonumber\\
I_4 = &  \int_\O \big(f(t,\bx,u_m(t,\bx))-f(t,\bx,u_n(t,\bx))\big)\phi_\eps(u_m(t,\bx)-u_n(t,\bx)\big) d\bx \nonumber
 \end{align*}
 For the rest of the proof, we use  $C$ to denote a constant which may change from one expression to another, but which is independent of $m,n$ and $\eps$. 
 
To estimate the integral $I_1$, we use the boundedness of $\phi_\eps$ in $L^\infty(\R)$  and the H\"older continuity of $A_i$   with respect to $s$, see \fer{p.a4}, together with \fer{p.a0}, which yield
\begin{align*}
I_1 \leq  C\nor{u_m}_{L^q}^{\alpha }\nor{J_m-J_n}_{L^\frac{q}{q-\alpha}} , 
\end{align*}
Since $\phi_\eps$ is non-decreasing and bounded in $L^\infty(\R)$, the term $I_2$ is non-positive due to the monotonicity of the approximants $A_i$ stated in  \fer{p.a1}, see Lemma~\ref{lemma:mon}. 
For $I_3$, we use the H\"older continuity of $A_i$ with respect to $i$, see \fer{p.a2},  and again the boundedness of $\phi_\eps$ in $L^\infty(\R)$, obtaining 
\begin{align*}
I_3 \leq C \nor{J_n}_{L^1}k(m,n),
\end{align*}
with $k(m,n)=(\abs{n-m}/\abs{mn})^\alpha$. 
Finally, since  $f(t,\bx,\cdot)$ is locally Lipschitz continuous we get, using $\abs{s\phi_\eps(s)} \leq 2\Phi_\eps(s)$, 
\begin{align*}
I_4 \leq L_f\int_\O \abs{u_m-u_n}\abs{\phi_\eps(u_m-u_n)} \leq 2L_f \int_\O \Phi_\eps(u_m-u_n).
\end{align*}
Using these estimates in \fer{me},  and that $u_m$ and $J_m$ are uniformly bounded in $L^1(Q_T)$ and ${L^\frac{q}{q-\alpha}}(B)$, respectively, we deduce 
\begin{align}
\p_t \int_\O \Phi_\eps(u_m-u_n) \leq C \big( \nor{J_m-J_n}_{L^\frac{q}{q-\alpha}}  + k(m,n)\big) + 2L_f \int_\O \Phi_\eps(u_m-u_n).
 \label{mon:1}
\end{align}
Since $J_m$ is a Cauchy sequence in $L^\frac{q}{q-\alpha}(B)$, for all  $\delta>0$ there exists $N>0$ such that 
\begin{align*}
\nor{J_m-J_n}_{L^\frac{q}{q-\alpha}} < \delta \qtext{and}\quad k(m,n)<\delta \qtext{for }m,n>N.
\end{align*}
Therefore, from \fer{mon:1} and Gronwall's lemma we get 
\begin{align}
\label{noreg1}
\int_\O \Phi_\eps(u_m(t,\cdot)-u_n(t,\cdot))  \leq & C \e^{2L_f T} \Big(  \int_\O \Phi_\eps(u_{0m}-u_{0n})   + \delta \Big).
\end{align}
Since $u_{0i}$ is a Cauchy sequence in $L^1(\O)$ and $\Phi_\eps \in C(\R)$ with $\Phi_\eps(0)=0$ and $\Phi_\eps\to \abs{\cdot}$ in $C(\R)$, we may 
redefine  $N>0$ to also have 
\begin{align*}
\int_\O \Phi_\eps(u_{0m}-u_{0n}) < \delta \text{ for }m,n>N,
\end{align*}
with $N$ independent of $\eps$. 
Using this bound and the theorem of dominated convergence in \fer{noreg1}, we deduce in the limit $\eps\to0$, 
\begin{align*}
\nor{u_m-u_n}_{L^1(Q_T)} \leq & C\delta,
\end{align*}
implying that $u_i$ is a Cauchy sequence in $L^1(Q_T)$. 
 Hence, there exists $u \in L^1(Q_T)$ such that $u_i \to u \qtext{strongly in } L^1(Q_T)$.
We then have, at least for a subsequence (not relabeled) that $u_i \to u$  a.e. in  $Q_T$.  
Since $u_i$ is uniformly bounded in $L^q(Q_T)$, the theorem of dominated convergence yields, if $q\in[1,\infty)$,
\begin{align*}
u_i \to u \qtext{strongly in } L^q(Q_T), 
\end{align*}
and, if $q=\infty$, 
\begin{align*}
u_i \to u \qtext{strongly in } L^r(Q_T), \text{ for any }r\in[1,\infty), \text{ with }u\in L^\infty(Q_T).
\end{align*}
In addition, we have directly from \fer{eq.eq.ap} and the uniform bounds of $J_i$ in $L^\frac{q}{q-\alpha}(B)$ and of ${u_i}$ in $L^q(Q_T)$  that 
$\nor{\p_t u_i}_{L^1}$ is uniformly bounded as well. Thus, at least for a subsequence (not relabeled), we have 
\begin{align*}
u_i \to u \qtext{weakly in } W^{1,1}(0,T;L^1(\O)).
\end{align*}
Therefore, replacing $u$ by $u_i$ in \fer{eq.eq.ap} and taking the limit $i\to\infty$, we find that the limit $u$ is a solution of \fer{eq.eq}-\fer{eq.id}. In addition, if $f(\cdot,\cdot,0)\geq 0$ then $u_i \geq 0$ and therefore we also have $u\geq 0$.

Finally, the stability result \fer{stability2} is easily deduced by modifying the argument employed in Theorem~\ref{th.existence-Lipschitz} to take into account the monotonicity of $A$. 
Using again $\phi(u)$, with $\phi(s)=\abs{s}^{q-1}s$, for $q\geq 1$, and $u=u_1-u_2$, as a test function in \fer{eq.stabi} we find, due to the monotonicity of $A$ and $\phi$, 
\begin{align}
\label{ineq.pp}
 \p_t \int \Phi(u(t,\bx))  \leq
\int_\O  \phi(u(t,\bx))  \big(f(t,\bx,u_1(t,\bx))- f(t,\bx,u_2(t,\bx))\big),
 \end{align}
 for $(t,\bx)\in Q_T$. 
 Then,  Gronwall's inequality implies
\begin{align*}
\nor{u}_{L^q(Q_T)} \leq \e^{L_fT} \nor{u_0}_{L^q(\O)}.
\end{align*}
Notice that if $f(t,\bx,\cdot)$ is non-increasing then we directly obtain from \fer{ineq.pp} and $\phi$ non-decreasing that 
$\nor{u}_{L^q(Q_T)} \leq \nor{u_0}_{L^q(\O)}$. 
$\Box$

\medskip 

\no\textbf{Proof of Lemma~\ref{lemma:mona}.}
Consider a mollifier $\rho_n\in C^\infty_c(\R)$ given by $\rho_n(s) = n \rho(ns)$ for $n\in \N$, for some even function $\rho \in C^\infty_c(\R)$ with $\rho\geq 0$ in $\R$ and  such that $\int_\R \rho(s)ds = 1$. Observe that since $\rho$ is of compact support,
\begin{align}
\label{mon.is}
I_\alpha = \int_\R \rho(s)\abs{s}^\alpha ds <\infty.
\end{align}
We then have that the sequence $A_n (t,\bx,\by,s)= \int_\R \rho_n(s-\sigma) A(t,\bx,\by,\sigma) d\sigma$ satisfies  \fer{p.a0} due to the even symmetry of $\rho$ and the odd symmetry of $A$. Since the variables $t,\bx,\by$ do not play any role in this proof, we omit them for clarity. We have 
\begin{align*}
A_n (s_1) - A_n (s_2) = & \int_\R \rho_n(\sigma) \big(A(s_1-\sigma) - A(s_2-\sigma)\big) d\sigma,
\end{align*}
from where the monotonicity of $A_n$ stated in  \fer{p.a1} is easily deduced from that of $A$.  
We also check from this identity that the H\"older continuity of $A_n$ with respect to $s$ stated in \fer{p.a4} holds with the same continuity constant than that of $A$, due to the normalization of $\rho_n$ assumed in \fer{mon.is}.
The H\"older continuity with respect to $n$ stated in  \fer{p.a2} is deduced as 
\begin{align*}
\left|A_{n_1} (s) - A_{n_2} (s)\right| \leq & \int_\R \rho(\xi) \left|(A\Big(s-\frac{\xi}{n_1}\Big) - A\Big(s-\frac{\xi}{n_2}\Big)\right| d\xi \\
& \leq C_H\left|\frac{1}{n_1}-\frac{1}{n_2}\right|^\alpha \int_\R \rho(\xi) \abs{\xi}^\alpha \leq I_\alpha \left|\frac{1}{n_1}-\frac{1}{n_2}\right|^\alpha,
\end{align*}
where $C_H$ is the H\"older continuity constant of $A$. 
Finally, for the convergence result \fer{p.a3}, we have   
\begin{align*}
\left|A_n (s) - A (s)\right| \leq & \int_\R \rho_n(\sigma) \left|A(s-\sigma) - A(s)\right| d\sigma \\
& \leq  C_H \int_\R \rho_n(\sigma) \abs{\sigma}^\alpha d\sigma
 = \frac{C_H}{n^\alpha} I_\alpha.
\end{align*}
$\Box$

\medskip

\no\textbf{Proof of Lemma~\ref{lemma:mon}.}
Using the identity \fer{prop.parts} of Lemma~\ref{lemma:intparts}, we get 
\begin{align*}
 \int_\O \phi\big(u(\bx)-v(\bx)\big)& (\cL(t) u(\bx)-\cL(t) v(\bx))d\bx =  -\frac{1}{2}\int_\O\int_\O \Big(J(\bx-\by) \\
& \times\big(A(t,\bx,\by,u(\by)-u(\bx)) - A(t,\bx,\by,v(\by)-v(\bx))\big) \\
 &\times \big(\phi(u(\by)-v(\by)) - \phi(u(\bx)-v(\bx))\big)\Big) d\by d\bx.
\end{align*}
Let  $\xi_0\in C^{0,\alpha}(\R)$ be non-decreasing and consider a sequence $\xi_\eps\in C_c^\infty(\R)$ such that $\xi'_\eps \geq0$  and $\xi_\eps\to\xi_0$ in $L^\infty(\R)$ as $\eps\to0$ (see the proof of \fer{p.a1} of Lemma~\ref{lemma:mona} for the construction of such sequence).
Let $I:\R\to\R$ be given by  
\begin{align*}
I(\eps)= \big( \phi(s_1-t_1) - \phi(s_2-t_2) \big) \big( \xi_\eps (s_1-s_2) - \xi_\eps(t_1-t_2) \big), 
 \end{align*}
for $t_i,s_i\in\R,~i=1,2$. The result of this lemma follows if we prove the pointwise bound $I(0)\geq0$.
We have 
\begin{align*}
 \xi_\eps (s_1-s_2) - \xi_\eps(t_1-t_2) = \xi_\eps'(\eta) (s_1-s_2-(t_1-t_2)),
\end{align*}
for some intermediate point $\eta$ between $s_1-s_2$ and $t_1-t_2$.  
 Since $\phi$ and $\xi_\eps$ are non-decreasing, we deduce
 \begin{align*}
I(\eps)=  \xi_\eps'(\eta)  \big( \phi(s_1-t_1) - \phi(s_2-t_2) \big)  (s_1-t_1- (s_2-t_2))\geq 0.
 \end{align*}
Taking the limit $\eps\to0$ we deduce $I(0)\geq0$. $\Box$

\end{document}